\documentclass[a4paper,14pt]{extarticle}
\usepackage[T1]{fontenc} 
\usepackage{ulem} 
\usepackage[utf8]{luainputenc} 
\usepackage{geometry} 
\usepackage[pdftex]{graphicx} 
\usepackage{amstext} 
\usepackage{amssymb}
\usepackage{amsmath} 
\usepackage[T1,T2A]{fontenc} 
\usepackage[utf8]{inputenc} 
\usepackage[english]{babel} 
\usepackage{setspace}
\usepackage{esint} 
\usepackage{comment}
\usepackage{babel} 
\usepackage{float} 
\usepackage{fullpage} 
\usepackage{tikz} 
\usepackage{indentfirst} 
\usepackage{mathtext}
\usepackage{ dsfont }
\usepackage{ amsfonts }
\usepackage{ mathrsfs }
\usepackage[utf8]{inputenc}

\newcommand\tab[1][1cm]{\hspace*{#1}}

\title{Convergence of some classes of random flights in Wasserstein distance. \footnote{The study has been funded by the Russian Academic Excellence Project '5-100}}
\author{
    A. R. Falaleev \footnote{Lomonosov Moscow State University, Moscow, Russian Federation}
    \and V. D. Konakov \footnote{Higher School of Economics, Pokrovsky boulevard 11, Moscow, Russian Federation. \newline \tab[0.5cm] VKonakov@hse.ru}
}
\date{}
\begin{document}

\maketitle

In this paper we consider a random walk of a particle in $\mathbb{R}^d$.
Convergence of different
transformations of trajectories of random flights with Poisson switching moments has been obtained in Davydov and Konakov [3], as well as diffusion approximation of the process has been built. The goal of this paper is to prove stronger convergence in terms of the Wasserstein distance. Three types of transformations are considered:
cases of exponential and super-exponential growth of a switching moment transformation function are quite simple, and the result follows from the fact that the limit processes belong
to the unit ball. In the case of the power growth the estimation is more complicated and follows from combinatorial reasoning and properties of the Wasserstein metric.


\section{Introduction}

Let $(\mathcal{X}, \textbf{d})$ be a Polish space, and  $p \in [1,\infty)$. Remind that the \newline Wasserstein space of order p is defined as 
$$
P_{p}(\mathcal{X}):=\left\{\mu \in P(\mathcal{X}) ; \quad \int_{\mathcal{X}} d\left(x_{0}, x\right)^{p} \mu(d x)<+\infty\right\}
$$
for some (and then for any) $x_0 \in \mathcal{X}$, where $P(\mathcal{X})$ is a set of all probability measures of $\mathcal{X}$ .

For any two probability measures $\mu$, $\nu$ on $\mathcal{X}$, the Wasserstein distance of order
p between $\mu$ and $\nu$ is defined by the formula

$$
W_{p}(\mu, \nu)=\left(\inf _{\pi \in \Pi(\mu, \nu)} \int_{\mathcal{X}} d(x, y)^{p} d \pi(x, y)\right)^{1 / p}
$$
$$
=\inf \left\{\left[\mathbb{E} \; \textbf{d}(X, Y)^{p}\right]^{\frac{1}{p}}, \quad \operatorname{law}(X)=\mu, \quad \operatorname{law}(Y)=\nu\right\}.
$$

It is well known that the Wasserstein distance $W_p$ metrizes the week convergence in $P_p(\mathcal{X})$ (Theorem 6.8 in [Vil06]). However, it is necessary to emphasize that the “weak convergence” in the sense of [Vil06] is stronger than the classical weak convergence (see the Definition 6.7 in [Vil06]). These two types of convergence coincide if the metric $d$ is bounded. But in general, they are different. The weak convergence established in [3] is in the classical sense. That is why the convergence in the Wasserstein metric should be proven and does not automatically follows from the weak convergence proved in [3]. 

Consider the random walk of a particle in $\mathbb{R}^d$ which is defined by two independent sequences of random variables $(T_k)$ and $(\varepsilon_k)$. The sequence $\varepsilon_k$ consists of independent random variables distributed on the unit sphere $S^{d-1}$  and defines the direction of motion of the particle. The sequence $T_k$, $ \forall k \mbox{ } T_k \ge 0$, $T_k \le T_{k+1}$ can be interpreted as moments when directions change. A particle starts from zero and moves in the direction $\varepsilon_1$ up to the moment $T_1$. It then
changes direction to $\varepsilon_2$ and moves on within the time interval $T_2 - T_1$, etc. The speed is
constant at all sites.  The position of the particle at the moment $t$ is denoted by $X (t)$.
In the article [3], the conditions under which the process $\{Y_T, T>0\}$
\[
Y_{T}(t) = \frac{1}{B(T)} X(tT), \mbox{ } t \in [0, 1]
\]
weakly converges in $C[0,1]$: $Y_T \Rightarrow Y$, $T \to \infty$ and $B_T \to \infty$ were \newline considered.

\par \medskip

The switching moments are assumed to form the Poisson process $\mathbb{T} = (T_k)$ in $\mathbb{R}_+$. In the homogeneous case, the process $X (t)$ is a random walk, because spacings $T_{k+1} - T_k$ are independent and $Y$ is a Wiener process.
In the case of non homogeneous Poisson process, the situation is getting complicated because of the increments $T_{k+1}-T_{k}$, which are not independent.

\par \medskip

Nevertheless,  the form of the limiting process was found and weak \newline convergence was proved for some Poisson switching moment transformation functions.
Let $T_k = f (\Gamma_k)$, where $(\Gamma_k)$ is a standard homogeneous Poisson process in $\mathbb{R}_{+}$ of intensity 1. In this case
\[
(\Gamma_k) = (\gamma_1 + \gamma_2 + ... + \gamma_k),
\]
where $(\gamma_k)$ are standart i.i.d. exponential random variables and $f(x)$ is a regular function with polynomial, exponential of super exponential growth. It is also assumed that $E\varepsilon_1 = 0$.
\par
We define a process
\[
Z_{n}(t) = Y_{T_n}(t).
\]
For $T = T_n$ the trajectories $\{Z_{n}(t)\mbox{ } t \in [0,1]\} $ are continuous broken lines with vertices at the points $\{(t_{n,k}, \frac{S_k}{B_n}), \mbox{ } k = 0, 1, ... , n\} $, where $t_{n,k} = \frac{T_k}{T_n},\mbox{ } T_0 = 0,\mbox{ } B_n = B(T_n),\mbox{ }
S_k = \sum_{i=1}^{k} \varepsilon_i (T_i - T_{i-1})$.
\\
The main result of the first part of [3] is following theorem:

\newtheorem{Tm}{Theorem}

\newtheorem{Df}{Definition}

\begin{Tm}
Under previous assumptions:

1) If the function $f$ has polynomial growth: $f(t) = t^{\alpha}, \alpha > 1/2$, we take $B(T) = T^{\frac{2\alpha-1}{2\alpha}}$. Then the process $Z_n$ converges weakly to $Y$, where $Y$ is a Gaussian process 
\[
Y(t) = \sqrt{2\alpha} \int_{0}^{t} s^\frac{\alpha-1}{2\alpha} dw(s),
\]
and $w$ is a process of Brownian motion for which the covariance matrix $w(1)$ coincides with the covariance matrix of $\varepsilon_1$
\par \medskip
2) If the function $f$ has exponential growth: $f(t) = e^{t\beta}, \beta > 0$, we take $B(T) = T$. Then the process $Z_n$ converges weakly to $Y$, where $Y$ is a continuous piecewise linear process with vertices at the points $(t_k, Y(t_k))$, \[
t_k = e^{-\beta\Gamma_{k-1}},\mbox{ } \Gamma_0 = 0,
\]
\[
Y(t_k) = \sum_{i=k}^{\infty} \varepsilon_i (e^{-\beta\Gamma_{i-1}} - e^{-\beta\Gamma_{i}}),\mbox{ } Y(0) = 0.
\]
\par \medskip
3) In the super-exponential case, suppose that $f$ is increasing continuous function such that
\[
\lim_{t \to \infty} \frac{f'(t)}{f(t)} = +\infty,
\]
We consider $B(T) = T$.
Then  $\frac{T_n}{T_{n+1}} \to 0$ in probability, and $Z_n \Rightarrow Y$, where the limit process degenerates:
\[
Y(t) = \varepsilon_1 t,\mbox{ } t \in [0, 1].
\]
\end{Tm}
\par
Recall the goal of this paper is to prove stronger convergence, namely convergence in Wasserstein distance. The values of constants may change from line to line.

\section{Main result}
In the rest of the paper we will consider $\mathcal{X} = C[0,1]$ and $\textbf{d}(f, g)=\sup_{x \in [0, 1]}(|f(x) - g(x)|)$. For a continuous random process $X(t), \; t \in [0,1]$, we denote a measure $\mu_X$ in $C[0,1]$ corresponding to this process.

\begin{Tm} Consider the Polish space $(\mathcal{X}, \textbf{d})$ and $p \in [1, \infty)$. Then

$$
W_p(\mu_{X_n}, \mu_Y) \rightarrow 0,
$$

where the process $X_n$ for the cases 1) - 3) is a polyline process with vertices at the points $(t_{n,k}, X_n(t_{n,k}))$. 

For the case 1)
$$
t_{n, k}=\left(\frac{\Gamma_{k}}{\Gamma_{n}}\right)^{\alpha}, \; X_{n}\left(t_{n, k}\right)=n^{\frac{1}{2}- \alpha} \sum_{i=1}^{k} \varepsilon_{i}\left(\Gamma_{i}^{\alpha}-\Gamma_{i-1}^{\alpha}\right), \; \Gamma_{0}^{\alpha}=0, \; k=1, \ldots, n.
$$
The limiting process $Y(t)$ is Gaussian
$$
Y(t) = \sqrt{2\alpha}\int_0^ts^{\frac{\alpha - 1}{2\alpha}}dw(s),
$$
where $w(s)$ is a process of Brownian motion for which the covariance matrix of $w(1)$ coincides with the covariance matrix of $\varepsilon_1$.
\newline
For the case 2)
\newline
$$
 t_{n, k}=e^{-\beta\left(\Gamma_{n}-\Gamma_{k}\right)}, \; X_{n}\left(t_{n, k}\right)=e^{-\beta \Gamma_{n}} \sum_{i=1}^{k} \varepsilon_{i}\left(e^{\beta \Gamma_{i}}-e^{\beta \Gamma_{i-1}}\right),    
$$
$$
\Gamma_{0}=0, \; k=1, \ldots, n.
$$

The limiting process $Y(t)$ is a continuous piecewise linear process with countable number of
vertices $\left(t_{k}, Y\left(t_{k}\right)\right),\; k=1,2, \ldots \:\;\; t_{k}=e^{-\beta \Gamma_{k-1}},\; \Gamma_{0}=0,$

$$
Y\left(t_{k}\right)=\sum_{i=k}^{\infty} \varepsilon_{i}\left(e^{-\beta \Gamma_{i-1}}-e^{-\beta \Gamma_{i}}\right).
$$
For the case 3)
$$
t_{n, k}=\frac{f\left(T_{k}\right)}{f\left(T_{n}\right)}, \; \;  X_{n}\left(t_{n, k}\right)=\frac{1}{f\left(T_{n}\right)} \sum_{i=1}^{k} \varepsilon_{i}\left(f\left(\Gamma_{i}\right)-f\left(\Gamma_{i-1}\right)\right),
$$
$$
\Gamma_{0}=0, k=1, \ldots, n.
$$
The limiting process $Y(t)$ degenerates:

$$
Y(t)=\varepsilon_{1} t, \;\;\; t \in[0,1].
$$

\end{Tm}
\section{Supporting definitions and results}
\begin{Df}{"Weak convergence in $P_p$"}

Let $(\mathcal{X}, d)$ be a Polish space, and $p \in
[1, \infty)$. Let $(\mu_k)k\in \mathbb{N}$ be a sequence of probability measures in $P_p(X)$ and let $\mu$ be another
element of $P_p(\mathcal{X})$. Then $\mu_k$ is said to "converge weakly in $P_p(X)$" if any from the following
equivalent properties is satisfied for some (and then any) $x_0 \in X$:
\par \medskip
i) $\mu_k \Rightarrow \mu$ when $k \to \infty$ and $$\int d(x_0, x)^pd\mu_{k}(x) \longrightarrow \int d(x_0, x)^pd\mu(x) \mbox{ } k \to \infty;$$
\par \medskip
ii) $\mu_k \Rightarrow \mu$ when $k \to \infty$ and $$\limsup _{k \rightarrow \infty} \int d\left(x_{0}, x\right)^{p} d \mu_{k}(x) \leq \int d\left(x_{0}, x\right)^{p} d \mu(x);$$
\par \medskip
iii) $\mu_k \Rightarrow \mu$ when $k \to \infty$ and
$$
\lim _{R \rightarrow \infty} \limsup _{k \rightarrow \infty} \int_{d\left(x_{0}, x\right) \geq R} d\left(x_{0}, x\right)^{p} d \mu_{k}(x)=0;
$$
\par \medskip
iv) For all continuous functions $\varphi$ with $|\varphi(x)| \leq C\left(1+d\left(x_{0}, x\right)^{p}\right), \\C\in \mathbb{R_+}$, one has 
$$
\int \varphi(x) d \mu_{k}(x) \longrightarrow \int \varphi(x) d \mu(x).
$$
\end{Df} 
\par \medskip
\begin{Tm}{($W_p$ metrizes $P_p(\mathcal{X})$, [1], Theorem 6.8).}
Let $(\mathcal{X}, d)$ be a Polish space, and $p \in
[1, \infty)$;  then
the Wasserstein distance $W_p$ metrizes the "weak convergence in $P_p(\mathcal{X})$". In other words, if
$(\mu_k)_{k\in\mathbb{N}}$ is a sequence of measures in $P_p(\mathcal{X})$ and $\mu$ is another measure in $P_p(\mathcal{X})$, then the statements
\begin{center}
$\mu_k$ "converges weakly in $P_p(\mathcal{X})$" to $\mu$
\end{center}
and
\begin{center}
$W_{p}\left(\mu_{k}, \mu\right) \longrightarrow 0$    
\end{center}
are equivalent.

\end{Tm}
To prove iii) in the Definition 1 we also need additional estimates. 
\begin{Tm}{(Doob's maximal inequality, [2])}.

If $X_k$ is a martingale or positive submartingale indexed by a finite set $ k \in (0,1, \ldots, N)$, then $ \forall p \geq 1$ and $\lambda>0$
\[
\lambda^{p} P\left[\sup _{0 \le k \le N}\left|X_{k}\right| \geq \lambda\right] \leq E\left[\left|X_{N}\right|^{p}\right],
\]
and for any $p>1$ 
\[ E\left[\left|X_{N}\right|^{p}\right] \leq E\left[\sup _{k}\left|X_{k}\right|^{p}\right] \leq\left(\frac{p}{p-1}\right)^{p} E\left[\left|X_{N}\right|^{p}\right].
\]
\end{Tm}

We will use the following estimates from [3].
\par \medskip
\noindent\textbf{Lemma 1}

Let $\alpha>0$ and $ m \geq 1$. Then $\forall x>0, h>0$
\[
{(x+h)^{\alpha}-x^{\alpha}=\sum_{k=1}^{m} a_{k} h^{k} x^{\alpha-k}+R(x, h)},
\]
where
\[{\quad a_{k}=\frac{\alpha(\alpha-1) \ldots(\alpha-k+1)}{k !}} ,
\]
\par \medskip
and
\[
{|R(x, h)| \leq\left|a_{m+1}\right| h^{m+1} \max \left\{x^{\alpha-(m+1)},(x+h)^{\alpha-(m+1)}\right\}}.
\]
\par \medskip
\noindent\textbf{Lemma 2}

For $\alpha \geq 0$ and $k \rightarrow \infty$ 
\[
{\quad\left(1+\frac{\alpha}{k}\right)^{k}=e^{\alpha}+O\left(\frac{1}{k}\right)}.
\]
\par \medskip
\noindent\textbf{Lemma 3}

Let $\Gamma$ denote the Gamma function. Then, when $k \rightarrow \infty$
\[
{\frac{\Gamma(k+\alpha)}{\Gamma(k)}=k^{\alpha}+O\left(k^{\alpha-1}\right)}.
\]
\par \medskip
\noindent\textbf{Lemma 4}

For any real $\beta$ at $k \rightarrow \infty$
\[
{E\left(\Gamma_{k}^{\beta}\right)=k^{\beta}+O\left(k^{\beta-1}\right)}.
\]
\par \medskip
\noindent\textbf{Lemma 5}

Let $\alpha \geq 0$. For $ k \rightarrow \infty$, the following relations hold:
\[
{\Gamma_{k+1}^{\alpha}-\Gamma_{k}^{\alpha}=\alpha \gamma_{k+1} \Gamma_{k}^{\alpha-1}+\rho_{k}},
\]
where $\left|\rho_{k}\right|=O\left(k^{\alpha-2}\right)$ in probability; 

\[
{\quad E\left|\Gamma_{k+1}^{\alpha}-\Gamma_{k}^{\alpha}\right|^{2}=2 \alpha^{2} k^{2 \alpha-2}+O\left(k^{2 \alpha-3}\right)}.
\]
From Lemma 5 we have

\par \medskip
\noindent\textbf{Corollary 1}
\[
{\sum_{1}^{n-1} E\left|\Gamma_{k+1}^{\alpha}-\Gamma_{k}^{\alpha}\right|^{2}=\frac{2 \alpha^{2}}{2 \alpha-1} n^{2 \alpha-1}+O\left(n^{2 \alpha-2}\right)}.
\]

\section{Proof of Theorem 2}
Let us consider three cases of the theorem.

\textit{Exponential growth case.}
Switching moment transformation function has the following form: $f(t) = e^{t\beta}, \beta > 0$ when $B(T) = T$, the process $Z_n$ converges weakly to $Y$, where $Y$ is a continuous piecewise linear process with vertices at the points $(t_k, Y(t_k))$, 
\[
t_k = e^{-\beta\Gamma_{k-1}},\mbox{ } \Gamma_0 = 0,
\]
\[
Y(t_k) = \sum_{i=k}^{\infty} \varepsilon_k (e^{-\beta\Gamma_{i-1}} - e^{-\beta\Gamma_{i}}),\mbox{ } Y(0) = 0.
\]
\par \medskip
For $T = T_n$ trajectories $\{Z_{n}(t)\mbox{ } t \in [0,1]\} $ are continuous broken lines with vertices $\{(t_{n,k}, \frac{S_k}{B_n}), \mbox{ } k = 0, 1, ... , n\} $, where $t_{n,k} = \frac{T_k}{T_n},\mbox{ } T_0 = 0,\mbox{ } B_n = B(T_n),\mbox{ } 
S_k = \sum_{i=1}^{k} \varepsilon_i (T_i - T_{i-1})$. 

Thus, the process takes the form of a polyline with nodes at points $(t_{n, k}, X_n(t_{n,k}))$:
\[
X_n(t_{n,k}) = \frac{1}{e^{\beta\Gamma_n}}\sum_{i=1}^{k} \varepsilon_i (e^{\beta\Gamma_{i}} - e^{\beta\Gamma_{i-1}}).
\]

The process $X_{n}(\cdot) \stackrel{\mathcal{L}}{=}Y_{n}(\cdot)$[(see [3], p.8)], where $Y_{n}(\cdot)$ is a polyline with the vertices $\left(\tau_{n, k}, Y_{n}\left(\tau_{n, k}\right)\right),\mbox{ }\left(\tau_{n, k}\right) \downarrow \;,\mbox{ } \tau_{n, 1}=1,\mbox{ } \tau_{n, k}=e^{-\beta\left(\gamma_{1}+\cdots+\gamma_{k-1}\right)},\mbox{ } k=2, \dots, n$,
$$
Y_{n}\left(\tau_{n, k}\right)=\sum_{i=k}^{n-1} \varepsilon_{i}\left(e^{-\beta\Gamma_{i-1}}-e^{-\beta\Gamma_{i}}\right)+\varepsilon_{n} e^{-\beta\Gamma_{n-1}},
$$
$Y_{n}(0)=0$ and $\Gamma_{0} =0$.

Since $Y_{n}\left (\tau_{n, k}\right)$ is a sum of non-negative terms multiplied by the random vector $\varepsilon_i$, $|\varepsilon_i | = 1$, then
$$
\max _{k = 1, \ldots , n}\left|Y_{n}(\tau_{n, k})\right| \leq \sum_{i=1}^{n-1} \left(e^{-\beta\Gamma_{i-1}}-e^{-\beta\Gamma_{i}}\right)+e^{-\beta\Gamma_{n-1}} = 1.
$$
Therefore for $R > 1$
$$
\mu_{n}(\textbf{d}(\textbf{0}, x) \geq R)=P\left(\max _{0 \leq t \leq 1}\left|Y_{n}(t)\right| \geq R\right) = 0.
$$ 
Then
$$
\lim _{R \rightarrow \infty} \overline{\lim _{n \rightarrow \infty}} \int_{d(0, x) \geq R} \textbf{d}^{p}(\textbf{0}, x) d \mu_{n}(x)=0.
$$
The convergence of $W_{p}\left(\mu_{n}, \mu\right) \rightarrow 0$ for any $p>1$ is proved.

\textit{Super exponential growth case.}

In this case, assuming $B_{n}=B\left(T_{n}\right)=T_{n} :$

$$
\max _{k = 1, \ldots ,  n}\left|X_{n}(t_{n, k})\right| \leq \sum_{k=1}^{n} \frac{T_{k}-T_{k-1}}{T_{n}}=\frac{T_{n}}{T_{n}} = 1.
$$

Therefore for $R > 1$
$$
\mu_{n}(\textbf{d}(\textbf{0}, x) \geq R)=P\left(\max _{0 \leq t \leq 1}\left|X_{n}(t)\right| \geq R\right) = 0.
$$ 
Then
$$
\lim _{R \rightarrow \infty} \overline{\lim _{n \rightarrow \infty}} \int_{d(0, x) \geq R} d^{p}(0, x) d \mu_{n}(x)=0.
$$
The convergence of the $W_{p}\left(\mu_{n}, \mu\right) \rightarrow 0$ for any $p > 1$ is proved.

\textit{Polynomial growth case.}

Note that if $\varepsilon_{j}$ is a uniformly distributed random variable on the unit ball in $R^d$, then $\left\langle\varepsilon_{i}, e_{j}\right\rangle$ is a one-dimensional random variable distributed symmetrically with respect to zero. Next, we use the fact that the odd moments of this variable are equal to 0. We have:
\[
{\mathrm{} \text  T_{k}=\Gamma_{k}^{\alpha},\mbox{ } \alpha>\frac{1}{2},\mbox{ } t_{n, k}=\frac{T_{k}}{T_{n}}=\left(\frac{\Gamma_{k}}{\Gamma_{n}}\right)^{\alpha},\mbox{ } B_{n}=n^{\alpha-1 / 2}, } 
\]
\[
 {\Gamma_0^\alpha= 0 , \;X_{n}\left(t_{n, k}\right)=\frac{1}{B_{n}} \sum_{i=1}^{k} \varepsilon_{i}\left(\Gamma_{i}^{\alpha}-\Gamma_{i-1}^{\alpha}\right)}.
\]
We have the following upper bound:

\[
P\left(\max _{k=1, \ldots, n}\left|\frac{1}{B_{n}} \sum_{i=1}^{k} \left\langle\varepsilon_{i}, e_{j}\right\rangle\left(\Gamma_{i}^{\alpha}-\Gamma_{i-1}^{\alpha}\right)\right| \geq 3 R\right) \leq
\]
\par \medskip

\[
\leq {P\left(\max _{k=1, \ldots, n}\left|\frac{1}{B_{n}} \sum_{i=1}^{k} \left\langle\varepsilon_{i}, e_{j}\right\rangle\left(\Gamma_{i}^{\alpha}-\Gamma_{i-1}^{\alpha}\right)-\frac{\alpha}{B_{n}} \sum_{i=1}^{k} \left\langle\varepsilon_{i}, e_{j}\right\rangle \gamma_{i} \Gamma_{i-1}^{\alpha-1}\right|>R\right)+} 
\]
\par \medskip
\[
{P\left(\max _{k=1, \ldots, n}\left|\frac{\alpha}{B_{n}} \sum_{i=1}^{k} \left\langle\varepsilon_{i}, e_{j}\right\rangle \gamma_{i} \Gamma_{i-1}^{\alpha-1}-\frac{\alpha}{B_{n}} \sum_{i=1}^{k} \left\langle\varepsilon_{i}, e_{j}\right\rangle \gamma_{i}(i-1)^{\alpha-1}\right|>R\right)+}
\]
\par \medskip
$$
+{P\left(\max _{k=1, \ldots, n}\left|\frac{\alpha}{B_{n}} \sum_{i=1}^{k} \left\langle\varepsilon_{i}, e_{j}\right\rangle \gamma_{i}(i-1)^{\alpha-1}\right|>R\right)=I+I I+I I I}.
$$
\newline
Estimation for I:

We use the Doob's Maximal inequality assuming $\lambda = B_nR$, $p = 2N$. $\mathfrak{M}_{n} = \sigma(\gamma_{1}, ... , \gamma_{n})$ is a filtration generated by $(\gamma_1 \cdots \gamma_n)$, then the process
$$
A_k^\alpha = \sum_{i=1}^{k} \left\langle\varepsilon_{i}, e_{j}\right\rangle\left(\Gamma_{i}^{\alpha}-\Gamma_{i-1}^{\alpha}-\alpha \gamma_{i} \Gamma_{i-1}^{\alpha-1}\right)
$$
becomes a conditional martingale. By the Doob's maximal inequality: 

$$
{I=E\left(P\left(\max _{k=1, \ldots, n}\left|A_k^\alpha\right|>B_{n} R\right) | \mathfrak{M}_{n}\right) \leq}
$$
\par \medskip
$$
\leq {\frac{1}{R^{2 N} n^{2 N \alpha-N}} E\left(A_k^\alpha\right)^{2 N}=}
$$
$$
={\frac{1}{R^{2 N} n^{2 N \alpha-N}} \sum_{k_{1}+\cdots+k_{n}=2 N} \frac{(2 N) !}{k_{1} ! \ldots k_{n} !} \prod_{i=1}^{n} E\left(\left\langle\varepsilon_{i}, e_{j}\right\rangle\right)^{k_{i}} E\left(\Gamma_{i}^{\alpha}-\Gamma_{i-1}^{\alpha}-\alpha \gamma_{i} \Gamma_{i-1}^{\alpha-1}\right)^{k_{i}}}.
$$

Note that if among $k_i$ there is at least one odd, then the corresponding summand in the sum is equal to zero due to the symmetry of the distribution. $\varepsilon_i$.

From (19) in [3] it follows that for $\frac{1}{2}<\alpha<2$
$$
\left|E\left(\Gamma_{i}^{\alpha}-\Gamma_{i-1}^{\alpha}-\alpha \gamma_{i} \Gamma_{i-1}^{\alpha-1}\right)^{k_{i}}\right| \leq E\left|\Gamma_{i}^{\alpha}-\Gamma_{i-1}^{\alpha}-\alpha \gamma_{i} \Gamma_{i-1}^{\alpha-1}\right|^{k_{i}} \leq
$$
$$
\leq {C(\alpha, N) E\left(\gamma_{i}\right)^{2 k_{i}} E \Gamma_{i-1}^{k_{i}(\alpha-2)} \leq C(\alpha, N)\left(2 k_{i}\right) ! i^{(\alpha-2) k_{i}}}.
$$
Therefore using that $\alpha < 2$
$$
\left|\prod_{i=1}^{n} E\left(\left\langle\varepsilon_{i}, e_{j}\right\rangle\right)^{k_{i}} E\left(\Gamma_{i}^{\alpha}-\Gamma_{i-1}^{\alpha}-\alpha \gamma_{i} \Gamma_{i-1}^{\alpha-1}\right)^{k_{i}}\right| \leq C(\alpha, N) \prod_{i=1}^{n} i^{k_{i}(\alpha-2)} \leq C(\alpha, N).
$$

For $\alpha \geq 2$ we use the Cauchy-Bunyakovsky-Schwarz inequality
$$
E\left(\gamma_{i}^{2 k_{i}} \Gamma_{i}^{(\alpha-2) k_{i}}\right) \leq \sqrt{E\left(\gamma_{i}\right)^{4 k_{i}}} \sqrt{E \Gamma_{i}^{(2 \alpha-4) k_{i}}} \leq C(\alpha, N) i^{(\alpha-2) k_{i}},
$$
and
$$
{\left|E\left(\Gamma_{i}^{\alpha}-\Gamma_{i-1}^{\alpha}-\alpha \gamma_{i} \Gamma_{i-1}^{\alpha-1}\right)^{k_{i}}\right| \leq E\left|\Gamma_{i}^{\alpha}-\Gamma_{i-1}^{\alpha}-\alpha \gamma_{i} \Gamma_{i-1}^{\alpha-1}\right|^{k_{i}} \leq} \\ {C(\alpha, N) i^{(\alpha-2) k_{i}}}.
$$
Therefore by Lemma 5,
$$
{\left|\frac{(2 N) !}{k_{1} ! \ldots k_{n} !} \prod_{i=1}^{n} E\left( \left\langle\varepsilon_{i}, e_{j}\right\rangle\right)^{k_{i}} E\left(\Gamma_{i}^{\alpha}-\Gamma_{i-1}^{\alpha}-\alpha \gamma_{i} \Gamma_{i-1}^{\alpha-1}\right)^{k_{i}}\right| \leq}
$$

$$
\leq {C(\alpha, N) \prod_{i=1}^{n} i^{(\alpha-2) k_{i}} \leq C(\alpha, N) \prod_{i=1}^{n} n^{(\alpha-2) k_{i}}=C(\alpha, N) n^{2 N \alpha-4 N}}.
$$
\newline
Now let us estimate the number of nonzero summands in the sum:
$
\sum_{k_{1}+\dots+k_{n}=2 N} \ldots
$
The bound $C(N)n^N$
for the number of summands is obtained by simple \newline combinatorial reasoning. 
We have:
\begin{equation}
I \leq \frac{C(\alpha, N)}{R^{2 N} n^{2 N \alpha-N}} n^{2 N \alpha-4 N} n^{N} < \frac{C(\alpha, N)}{R^{2 N}}.
\end{equation}

This bound is sufficient to verify (iii) of the Definition 1. 

Estimation for II:

Similarly, we use the independence of $\gamma_i$ and $\Gamma_{i-1}^\alpha$ to obtain:
$$
I I \leq \frac{\alpha^{2 N}}{R^{2 N} n^{2 N \alpha-N}} \sum_{k_{1}+\cdots+k_{n}=2 N} D_{\overline{k}},
$$
where
$$
D_{\overline{k}} =\frac{(2 N) !}{k_{1} ! \ldots k_{n} !} \prod_{i=2}^{n} E\left(\left\langle\varepsilon_{i}, e_{j}\right\rangle\right)^{k_{i}} E\left(\gamma_{i}\right)^{k_{i}} E\left(\Gamma_{i-1}^{\alpha-1}-(i-1)^{\alpha-1}\right)^{k_{i}},
$$
$\overline{k} = (k_1, \cdots, k_n).$ Estimate the expectation:
$$
E\left(\Gamma_{i-1}^{\alpha-1}-(i-1)^{\alpha-1}\right)^{k_i}=\sum_{m=0}^{k_i}C_{k_i}^{m}(-1)^{k_i-m}(i-1)^{(k_i-m)(\alpha-1)}  E \Gamma_{i-1}^{(\alpha-1) m}= 
$$

$$
={\sum_{m=0}^{k_i}(-1)^{k_i-m}(i-1)^{(k_i-m)(\alpha-1)} C_{k_i}^{m}\left[(i-1)^{(\alpha-1) m}+O_{m}\left((i-1)^{(\alpha-1) m-1}\right)=\right.}
$$
$$
= {(i-1)^{k_i(\alpha-1)} [\sum_{m=0}^{k_i} C_{k_i}^{m}(-1)^{k_i-m}+ \sum_{m=0}^{k_i} C_{k_i}^{m}(-1)^{k_i-m} O_{m}\left((i-1)^{-1}\right)}].
$$
Now we have:
$$
{\left|\prod_{i=2}^{n} E\left(\left\langle\varepsilon_{i}, e_{j}\right\rangle\right)^{k_{i}} E\left(\gamma_{i}\right)^{k_{i}} E\left(\Gamma_{i-1}^{\alpha-1}-(i-1)^{\alpha-1}\right)^{k_{i}}\right| \leq C(N, \alpha) \prod_{i=1}^{n} i^{k_{i}(\alpha-1)} \leq} 
$$
$$
\leq {C(N, \alpha) \prod_{i=1}^{n} n^{k_{i}(\alpha-1)}=C(N, \alpha) n^{2 N(\alpha-1)}}.
$$
\\
We get:
\begin{equation}
{I I \leq C(N, \alpha) \frac{\alpha^{2 N}}{R^{2 N} n^{2 N \alpha-N}} n^{N} n^{2 N(\alpha-1)-1} \le C(N, \alpha) \frac{\alpha^{2 N}}{R^{2 N}}}.
\end{equation}

Now we have:
$$
I I I \leq \frac{\alpha^{2 N}}{R^{2 N} n^{2 N \alpha-N}} \sum_{k_{1}+\cdots+k_{n}=2 N} \frac{(2 N) !}{k_{1} ! \ldots k_{n} !} \prod_{i=2}^{n} E\left(\left\langle\varepsilon_{i}, e_{j}\right\rangle\right)^{k_{i}} E\left(\gamma_{i}\right)^{k_{i}}(i-1)^{k_{i}(\alpha-1)},
$$
where
$$
\left|\prod_{i=2}^{n} E\left(\left\langle\varepsilon_{i}, e_{j}\right\rangle\right)^{k_{i}} E\left(\gamma_{i}\right)^{k_{i}}(i-1)^{k_{i}(\alpha-1)}\right| \leq C(N, \alpha) n^{2 N(\alpha-1)}.
$$
We finally get:
\begin{equation}
I I I \leq C(N, \alpha) \frac{1}{R^{2 N} n^{2 N \alpha-N}} n^{N} n^{2 N(\alpha-1)}=\frac{C(N, \alpha)}{R^{2 N}}.
\end{equation}

The multidimensional case reduces to the 1-dimensional bounds in the following way:

$$
P\left(\max _{k=1,2, \ldots, n}\left|\sum_{i=1}^{k} \varepsilon_{i}\left(\Gamma_{i}^{\alpha}-\Gamma_{i-1}^{\alpha}\right)\right| \geq B_{n} R\right)=
$$
$$
= P\left(\max _{k=1,2, \ldots, n}\left|\sum_{j=1}^{d} \sum_{i=1}^{k}\left\langle\varepsilon_{i}, e_{j}\right\rangle e_{j}\left(\Gamma_{i}^{\alpha}-\Gamma_{i-1}^{\alpha}\right)\right| \geq B_{n} R\right) \leq
$$
$$
\leq P\left(\exists j^{*}, 1 \leq j^{*} \leq d,  \max _{k=1,2, \ldots, n}\left|\sum_{i=1}^{k}\left\langle\varepsilon_{i}, e_{j^{*}}\right\rangle e_{j^{*}}\left(\Gamma_{i}^{\alpha}-\Gamma_{i-1}^{\alpha}\right)\right| \geq \frac{B_{n} R}{d}\right) \leq
$$

\begin{equation}
    \le \sum_{j=1}^{d} P\left(\max _{k=1,2, \ldots, n}\sum_{i=1}^{k}\left|\left\langle\varepsilon_{i}, e_{j}\right\rangle\right|\left(\Gamma_{i}^{\alpha}-\Gamma_{i-1}^{\alpha}\right) \geq \frac{B_{n} R}{d}\right).
\end{equation}

Polynomial case estimation:

We check now the conditions of Theorem 6.8 from [1] using previously obtained
bounds (1) – (3) and (4):

$$
\int_{d(0, x)>R} \textbf{d}^{p}(\textbf{0}, x) d \mu_{n}(x) = \sum_{i=0}^{\infty} \int_{(i+1) R \leq \textbf{d}(\textbf{0}, x)<(i+2) R} \textbf{d}^{p}(0, x) d \mu_{n}(x) \leq
$$
$$
     \leq R^{p} \sum_{i=0}^{\infty}(i+2)^{p} \cdot \mu_{n}(\textbf{d}(\textbf{0}, x) \geq(i+1) R) = 
     $$
     $$
     =  R^{p} \sum_{i=0}^{\infty}(i+2)^{p} \cdot P\left(\max _{0 \leq t \leq 1}\left|X_{n}(t)\right| \geq (i+1)R\right) =
$$
$$
     = R^{p} \sum_{i=0}^{\infty}(i+2)^{p} \cdot P\left(\max _{k=1, \ldots, n}\left|\frac{1}{B_{n}} \sum_{i=1}^{k} \varepsilon_{i}\left(\Gamma_{i}^{\alpha}-\Gamma_{i-1}^{\alpha}\right)\right| \geq (i+1)R\right)\leq
$$
$$
\leq  R^{p} \sum_{i=0}^{\infty}(i+2)^{p} \cdot \sum_{j=1}^{d} P\left(\max _{k=1,2, \ldots, n}\sum_{i=1}^{k}\left|\left\langle\varepsilon_{i}, e_{j}\right\rangle\right|\left(\Gamma_{i}^{\alpha}-\Gamma_{i-1}^{\alpha}\right) \geq \frac{B_{n}(i+1) R}{d}\right) \leq
$$
$$
\le \frac{C(N, \alpha, d)}{R^{2N-p}}.
$$

Thus the condition iii) from the Definition 1 is satisfied for $\alpha > \frac{1}{2}$:
$$
\lim _{R \rightarrow \infty} \overline{\lim _{n \rightarrow \infty}} \int_{d(0, x) \geq R} d^{p}(0, x) d \mu_{n}(x)=0.
$$
This completes the proof of convergence in Wasserstein distance for $\alpha > \frac{1}{2}$ and $p \in [1, \infty)$.

The authors are grateful to Yu. Davydov for reading the article and \newline providing useful comments.

\section{References}
\renewcommand\refname{}

\end{document}